# Series solutions of Bessel-type differential equation in terms of orthogonal polynomials and physical applications

A. D. Alhaidari[a] and H. Bahlouli[b]

[a] *Saudi Center for Theoretical Physics, P.O. Box 32741, Jeddah 21438, Saudi Arabia*

[b] *Physics Department, King Fahd University of Petroleum & Minerals, Dhahran 31261, Saudi Arabia*

**Abstract:** We obtain a class of exact solutions of a Bessel-type differential equation, which is a six-parameter linear ordinary differential equation of the second order with irregular (essential) singularity at the origin. The solutions are obtained using the Tridiagonal Representation Approach (TRA) as bounded series of square integrable functions written in terms of the Bessel polynomial on the real line. The expansion coefficients of the series are orthogonal polynomials in the equation parameters space. We use our findings to obtain solutions of the Schrödinger equation for some novel potential functions.

**MSC**: 34-xx, 33C45, 33D45, 34B30, 81Qxx

**Keywords**: ordinary differential equations, Bessel polynomial, tridiagonal representation approach, recursion relation, hypergeometric polynomials, novel potentials

*Dedicated to the memory of our friend, colleague and collaborator,*

*the late Mohammed S. Abdelmonem.*

## 1. Introduction

It is common knowledge that differential equations play an important role in describing a large number of physical phenomena. Hence, it comes as no surprise that theoretical physicists spend much of their efforts in designing tools that enable them to find solutions of such differential equations. On the other hand, singularities of the differential equations and their types play a pivotal role in the nature of physical phenomena and controlling its richness. In the present work, we use a special approach to express the exact solution of a six-parameter Bessel-type differential equation as bounded series (finite or infinite) of a square integrable basis set.

The "Tridiagonal Representation Approach" (TRA) is an algebraic method that was established recently as a physics tool to solve quantum mechanical problems [1]. It was inspired by the J-matrix method of scattering [2-4]. The mathematical foundation of the method was established by Ismail and Koelink [4,5]. Additionally, the algebraic theory of tridiagonalization was further developed in the work by V. X. Genest *et. al* [6]. One of the advantages of the TRA as a physics tool is that its solution space is larger than that of the conventional methods which usually provide analytical treatment only for a limited number of well-known exactly solvable problems (e.g., the oscillator, Coulomb, Pöschl-Teller, Morse, etc.) [7,8]. As a mathematical tool, the TRA has recently been used to solve linear ordinary differential equations of the second order [9-11] of which the Schrödinger equation is an example. The solutions are written as bounded series (finite or infinite) of square integrable functions provided that the expansion coefficients of the series satisfy three-term recursion relations. These relations are then solved in terms of orthogonal polynomials on the real line with arguments and parameters related to the differential equation parameters.

In this work, we extend the TRA treatment to the solution of the following six-parameter differential equation

$$\left[ x^2 \frac{d^2}{dx^2} + (ax+b)\frac{d}{dx} + A_+ x + \frac{A_-}{x} + \frac{A_1}{x^2} - A_0 \right] y(x) = 0, \tag{1}$$

where $x \geq 0$ and $\{a, b, A_\pm, A_1, A_0\}$ are real parameters. This equation has an essential (irregular) singularity at $x = 0$ and if $A_+ \neq 0$ then it has another at infinity. The differential equation of the Bessel polynomial $J_n^\mu(x)$ is a special case of (1) with $A_\pm = A_1 = 0$, $a = 2(\mu+1)$, $b = 1$ and $A_0 = n(n+2\mu+1)$ (see pages 245 in [12]). In this work, we search for TRA solutions of Eq. (1) that could be written as follows

$$y(x) = \sum_n f_n \phi_n(x), \tag{2}$$

where $\{\phi_n(x)\}$ is a complete set of square integrable functions and $\{f_n\}$ are the expansion coefficients. Applying Fuch's theorem to Eq. (1) with its singularity structure as noted above, allows us to use Frobenius method and propose the following elements of the basis set

$$\phi_n(x) = G_n x^\alpha e^{-\beta/x} J_n^\mu(x), \tag{3}$$

where $J_n^\mu(x)$ is the Bessel polynomial on the real line whose properties are given in Appendix A and $G_n$ is a conveniently chosen normalization constant. The parameters $\{\alpha, \beta, \mu\}$ will be related to the differential equation parameters by the TRA constraints as will be shown below. However, $\mu$ should always be a negative real number. In section 2, we substitute the series (2) with the basis elements (3) into the differential equation (1) and identify the scenarios where three-term recursion relations for the expansion coefficients $\{f_n\}$ are produced. In sections 3 to 6, we obtain the recursion coefficients for each scenario and identify the orthogonal polynomials that satisfy the corresponding three-term recursion relation. As such, the solution of (1) is fully determined. Finally, in section 7, we provide a physical application of our findings where we transform the Schrödinger equation into Eq. (1) and hence identify the corresponding potential function, energy and wavefunction.

## 2. The TRA setup

If we write Eq. (1) as $\mathcal{D} y(x) = 0$, then the action of the second order differential operator $\mathcal{D}$ on the basis elements (3) could be written as follows

$$\mathcal{D}\phi_n(x) = G_n x^\alpha e^{-\beta/x} \left\{ x^2 \frac{d^2}{dx^2} + \left[ x(2\alpha + a) + 2\beta + b \right] \frac{d}{dx} \right.$$
$$\left. + A_+ x + \frac{1}{x}\left[ A_- + \alpha b + \beta(2\alpha + a - 2) \right] + \frac{A_1 + \beta(\beta + b)}{x^2} - A_0 + \alpha(\alpha + a - 1) \right\} J_n^\mu(x) \tag{4}$$

Using the differential equation of the Bessel polynomial (A4) in Appendix A, we can write this equation as

$$\mathcal{D}\phi_n(x) = G_n x^\alpha e^{-\beta/x} \left\{ \left[ 2x\left(\alpha - \mu - 1 + \tfrac{a}{2}\right) + 2\beta + b - 1 \right] \frac{d}{dx} + n(n + 2\mu + 1) \right.$$
$$\left. + A_+ x + \frac{1}{x}\left[ A_- + \alpha b + \beta(2\alpha + a - 2) \right] + \frac{A_1 + \beta(\beta + b)}{x^2} - A_0 + \alpha(\alpha + a - 1) \right\} J_n^\mu(x) \tag{5}$$

Now, the TRA requires that $\mathcal{D}\phi_n(x)$ be of the following form [1]

$$\mathcal{D}\phi_n(x) = \omega(x)\left[ u_n \phi_n(x) + s_{n-1} \phi_{n-1}(x) + t_n \phi_{n+1}(x) \right], \tag{6}$$

where $\omega(x)$ is non-zero entire function on the positive real line and the coefficients $\{u_n, s_n, t_n\}$ are $x$-independent. Substituting (2) in (1), written as $\mathcal{D}y(x) = 0$, and using (6) leads to the following sought-after three-term recursion relation for the expansion coefficients $\{f_n\}$

$$u_n Q_n + t_{n-1} Q_{n-1} + s_n Q_{n+1} = 0, \tag{7a}$$

where we have written $f_n = f_0 Q_n$, which makes $Q_0 = 1$. On the other hand, if we define $P_n = \left( \prod_{m=0}^{n-1} s_m / t_m \right) Q_n$ with $P_0 = Q_0 = 1$, then we can rewrite (7a) in the following alternative but useful form

$$u_n P_n + s_{n-1} P_{n-1} + t_n P_{n+1} = 0, \tag{7b}$$

which is obtained from (7a) by the exchange $s_n \leftrightarrow t_n$. Such a transformation turns out to be very useful in many circumstances. If we write $u_n$ as $a_n - zc_n$ such that $\{a_n, c_n, s_n, t_n\}$ are independent of $z$ then $P_n(z)$ will be a polynomial in $z$ of degree $n$. Throughout this work, we will find out that $c_n$ is always a constant independent of $n$. Moreover, if $s_n t_n > 0$ for all $n$ then the recursion is said to be definite (i.e., $s_n$ and $t_n$ have the same sign) and, in accordance with Favard's theorem [13], $\{Q_n(z)\}$ and $\{P_n(z)\}$ are orthogonal polynomials on the real line with $f_0^2(z)$ as their positive definite weight function [1,14,15].

Investigation shows that for a fixed lucid choice of $\omega(x)$ in front of the curly brackets in (5), there are only two actions that will turn (5) into (6). One, is the recursion relation of the Bessel polynomial (A2), which allows only linear function of $x$ inside the curly brackets to multiply $J_n^\mu(x)$. The other is the differential relation (A8), which allows only the differential operator $x^2 \frac{d}{dx}$ inside the curly brackets to act on $J_n^\mu(x)$. All other terms must be eliminated by imposing essential parameter relations or parameter constraints. Following this procedure, we arrive at only three possibilities each of which is endowed with a solution class:

(1) $\alpha = \mu+1-\frac{a}{2}$ and $2\beta = 1-b$: with $\omega(x) = x^{-k}$ where $k = 0,1,2,$ or $3$. (8a)

(2) $\alpha \neq \mu+1-\frac{a}{2}$ and $2\beta = 1-b$: with $\omega(x) = x^{-1}$. (8b)

(3) $\alpha = \mu+1-\frac{a}{2}$ and $2\beta \neq 1-b$: with $\omega(x) = x^{-2}$. (8c)

In the following three sections, we treat these cases separately. The corresponding three-term recursion relation is obtained and solved in terms of orthogonal polynomials.

## 3. TRA solution class corresponding to (8a)

In this case $\omega(x) = x^{-k}$ for $k = 0,1,2,$ or $3$ and Eq. (5) reduces to the following:

$$\mathcal{D}\phi_n(x) = G_n x^{\alpha-k} e^{-\beta/x} \left\{ x^{k-2}\left[A_1 + \tfrac{1}{4}(1-b^2)\right] + x^{k-1}\left[A_- + b\left(1-\tfrac{a}{2}\right) + \mu\right] + \right. $$
$$\left. + x^k\left[-A_0 - \tfrac{1}{4}(a-1)^2 + \left(n+\mu+\tfrac{1}{2}\right)^2\right] + A_+ x^{k+1}\right\} J_n^\mu(x) \quad (9)$$

In the following four subsections, we consider the case for each $k$ separately.

### 3.1: $k = 0$

In addition to the parameter constraints given by (8a), we must impose two additional constraints to eliminate the two non-tridiagonal terms inside the curly brackets in (9) which are proportional to $x^{-1}$ and $x^{-2}$. The result is the following parameter relations:

$$\alpha = (b-1)\left(\tfrac{a}{2}-1\right) - A_-,\ 2\beta = 1-b,\ \mu = b\left(\tfrac{a}{2}-1\right) - A_-,\ \text{and}\ b^2 = 1+4A_1. \quad (10)$$

The first three relations express the basis parameters in terms of the differential equation parameters. On the other hand, the last is a restriction on the differential equation parameters to obtain a solution. It relates $b$ to $A_1$ and thus reduces the number of free parameters of the differential equation from six to five while also imposing the reality constraint that $A_1 \geq -\tfrac{1}{4}$. With all of these constraints that result from the tridiagonal representation requirement, equation (5) transforms into the following:

$$\mathcal{D}\phi_n(x) = G_n x^\alpha e^{-\beta/x}\left[A_+ x + \left(n+\mu+\tfrac{1}{2}\right)^2 - v^2\right] J_n^\mu(x), \quad (11)$$

where $v^2 = A_0 + \tfrac{1}{4}(a-1)^2$. Using the recursion relation of the Bessel polynomial (A2) for the term $A_+ x J_n^\mu(x)$, this equation becomes

$$\mathcal{D}\phi_n(x) = G_n x^\alpha e^{-\beta/x} \frac{A_+}{2} \left\{ \left[ \frac{-\mu}{(n+\mu)(n+\mu+1)} + \frac{2}{A_+} \left( \left(n+\mu+\tfrac{1}{2}\right)^2 - \nu^2 \right) \right] J_n^\mu(x) \right.$$
$$\left. - \frac{n}{(n+\mu)(2n+2\mu+1)} J_{n-1}^\mu(x) + \frac{n+2\mu+1}{(n+\mu+1)(2n+2\mu+1)} J_{n+1}^\mu(x) \right\} \quad (12)$$

If we choose a basis normalization such that $G_n = 1$ and then compare (12) to (6), we obtain $\omega(x) = A_+/4$ and the following recursion coefficients:

$$u_n = \frac{-2\mu}{(n+\mu)(n+\mu+1)} + \frac{4}{A_+}\left[\left(n+\mu+\tfrac{1}{2}\right)^2 - \nu^2\right], \quad (13a)$$

$$s_n = -\frac{n+1}{(n+\mu+1)\left(n+\mu+\tfrac{3}{2}\right)}, \qquad t_n = \frac{n+2\mu+1}{(n+\mu+1)\left(n+\mu+\tfrac{1}{2}\right)}. \quad (13b)$$

Comparing the recursion relation (7b) constructed from these coefficients to (B1) in Appendix B, we can identify $P_n$ with the newly introduced polynomial $B_n^\mu(z;\gamma)$ [16] having the following argument and parameter

$$z = 4\nu^2/A_+ \text{ and } \gamma = 4/A_+. \quad (14)$$

Finally, the un-normalized solution (without the overall factor $F_0$) of the differential equation (1) for this class is written as follows

$$y(x) = \sum_n C_n B_n^\mu(z;\gamma) \phi_n(x), \quad (15)$$

where $C_n = \left(\prod_{m=0}^{n-1} t_m/s_m\right) = \left[\dfrac{n+\mu+\tfrac{1}{2}}{\mu+\tfrac{1}{2}}\right]\dfrac{(-n-2\mu)_n}{n!}$. Note that the normalized solution must be multiplied by $f_0 = \sqrt{\rho(z)}$ with $\rho(z)$ being the weight function for the orthogonal polynomial $B_n^\mu(z;\gamma)$, which is not known yet.

**3.2: $k = 1$**

Here too, additional parameter constraints, on top of those already given by (8a), should be imposed. These are needed to eliminate the two non-tridiagonal terms inside the curly brackets in (9) which are proportional to $x^2$ and $x^{-1}$. The result is the following parameter relations:

$$\alpha = \mu + 1 - \tfrac{a}{2}, \qquad 2\beta = 1 - b, \qquad b^2 = 1 + 4A_1, \qquad A_+ = 0 \quad (16a)$$

The first two relates the basis parameters to the differential equation parameters leaving the basis parameter $\mu$ free to be determined later in section 7 by physical requirements. On the other hand, the last two relations are restrictions on the differential equation parameters to obtain a solution. The third relates $b$ to $A_1$ while also imposing the reality constraint that $A_1 \geq -\tfrac{1}{4}$. The last relation eliminates the $A_+$ term from the differential equation making the singularity at infinity regular. Therefore, the number of free parameters of the differential equation are reduced from six to four. With all of these tridiagonal representation constraints, equation (5) becomes:

$$\mathcal{D}\phi_n(x) = G_n x^{\alpha-1} e^{-\beta/x} \left\{ x\left[\left(n+\mu+\tfrac{1}{2}\right)^2 - v^2\right] + \xi + \mu \right\} J_n^\mu(x), \tag{17}$$

where $\xi = A_- + b\left(1 - \tfrac{a}{2}\right)$. Using the recursion relation of the Bessel polynomial (A2) for the first term, this equation becomes

$$\mathcal{D}\phi_n(x) = \tfrac{1}{2} G_n x^{\alpha-1} e^{-\beta/x} \left\{ \left[ 2(\xi+\mu) - \mu \frac{\left(n+\mu+\tfrac{1}{2}\right)^2 - v^2}{(n+\mu)(n+\mu+1)} \right] J_n^\mu(x) \right.$$
$$\left. + \left[\left(n+\mu+\tfrac{1}{2}\right)^2 - v^2\right]\left[ -\frac{n J_{n-1}^\mu(x)}{(n+\mu)(2n+2\mu+1)} + \frac{(n+2\mu+1) J_{n+1}^\mu(x)}{(n+\mu+1)(2n+2\mu+1)} \right] \right\} \tag{18}$$

If we choose a basis normalization such that $G_n = 1$ and then compare (18) to (6), we obtain $\omega(x) = 1/4x$ and the following recursion coefficients:

$$u_n = 4(\xi+\mu) - 2\mu \frac{\left(n+\mu+\tfrac{1}{2}\right)^2 - v^2}{(n+\mu)(n+\mu+1)}, \tag{19a}$$

$$s_n = -\frac{n+1}{(n+\mu+1)\left(n+\mu+\tfrac{3}{2}\right)}\left[\left(n+\mu+\tfrac{3}{2}\right)^2 - v^2\right], \tag{19b}$$

$$t_n = \frac{n+2\mu+1}{(n+\mu+1)\left(n+\mu+\tfrac{1}{2}\right)}\left[\left(n+\mu+\tfrac{1}{2}\right)^2 - v^2\right]. \tag{19c}$$

Comparing the recursion relation (7b) constructed from these coefficients to (B7) in Appendix B, we conclude that $P_n$ is a special case of the discrete Hahn polynomial $Q_n^N(k;p,q)$ with

$$p = \mu - v - \tfrac{1}{2}, \qquad q = \mu + v + \tfrac{1}{2} = -N, \qquad k = -(\mu+\xi) \tag{20}$$

Therefore, the un-normalized solution of the differential equation (1) for this class reads as follows

$$y_k(x) = \sum_n C_n Q_n^N(k;p,q) \phi_n(x), \tag{21}$$

where $C_n = \left(\prod_{m=0}^{n-1} t_m/s_m\right) = \frac{\left(\mu+\tfrac{1}{2}\right)^2 - v^2}{\mu+\tfrac{1}{2}} \left[\frac{n+\mu+\tfrac{1}{2}}{\left(n+\mu+\tfrac{1}{2}\right)^2 - v^2} \frac{(-n-2\mu)_n}{n!}\right].$

On the other hand, if we compare the recursion relation (7b) constructed from the coefficients (19) to (B9) but for $(-1)^n H_n^p(z;q,c,d)$, then we get $p = -(\mu+\xi)$, $q = 1-(\mu+\xi)$, $c = 2\mu+\xi+\tfrac{1}{2}+v$, $d = 2\mu+\xi+\tfrac{1}{2}-v$ and $z = 0$.

**3.3: $k = 2$**

Here too, additional parameter constraints on top of those already given by (8a) should be imposed. These are needed to eliminate the two non-tridiagonal terms inside the curly brackets in (9) which are proportional to $x^2$ and $x^3$. The result is the following parameter relations:

$$\alpha = \mu + 1 - \tfrac{a}{2}, \qquad 2\beta = 1 - b, \tag{22a}$$

$$\left(n + \mu + \tfrac{1}{2}\right)^2 = A_0 + \tfrac{1}{4}(a-1)^2, \qquad A_+ = 0. \tag{22b}$$

The first relation in (22b) dictates that either the equation parameters depend on $n$ (e.g., $A_0$ is quadratic in $n$ or $a$ is linear in $n$) or $\mu$ should be written as $\mu(n) = -n - \nu - \tfrac{1}{2}$ where $\nu$ is the constant $\sqrt{A_0 + \tfrac{1}{4}(a-1)^2}$. The first choice is not acceptable because if the equation parameter(s) depend on $n$, then the solution should be $y(x) \propto \phi_n(x)$ and we must impose the diagonal rather than the tridiagonal structure, which will also requires that $A_1 + \beta(\beta+b) = 0$. Consequently, we end up with the trivial result that the equation becomes that of the Bessel polynomial $J_n^\mu(x)$. This is also evident by substituting the parameters from (22) in addition to $A_1 + \beta(\beta+b) = 0$ into Eq. (4). Therefore, the parameter relations (22) must be rewritten as

$$-\alpha(n) = n + \nu - \tfrac{1}{2} + \tfrac{a}{2}, \qquad -\mu(n) = n + \nu + \tfrac{1}{2}, \qquad 2\beta = 1-b, \qquad A_+ = 0, \tag{23}$$

where $\nu = \sqrt{A_0 + \tfrac{1}{4}(a-1)^2}$, which dictates that $4A_0 \geq -(a-1)^2$. Due to the dependence of the parameter $\mu$ on $n$, the form of the Bessel polynomial and its properties change. If we write $J_n^{\mu(n)}(x) := \bar{J}_n^\nu(x)$, then $\bar{J}_n^\nu(x) = n!(-x)^n L_n^{2\nu}(1/x)$ and its properties are given in Appendix A. Using the parameter relations (23) turns equation (9) into the following:

$$\mathcal{D}\phi_n(x) = G_n x^{\alpha-2} e^{-\beta/x} \left\{ \left[\xi - \left(n+\nu+\tfrac{1}{2}\right)\right]x + \left[A_1 + \tfrac{1}{4}(1-b^2)\right] \right\} \bar{J}_n^\nu(x). \tag{24}$$

Now, the properties of the polynomial $\bar{J}_n^\nu(x)$ given by relations (A11) and (A12) show that this case will not produce a tridiagonal structure of the form required by Eq, (6). On the other hand, since the n-dependent expression in $\phi_n(x)$ for this case is $G_n x^{-n} \bar{J}_n^\nu(x)$, which is equal to $g_n L_n^{2\nu}(1/x)$ with $g_n = (-1)^n n! G_n$, then the recursion relation of the Laguerre polynomial dictates that we must multiply $\phi_n(x)$ by $1/x$ not by $x$ to obtain the tridiagonal structure. Therefore, a proper treatment of this case should be conducted in another basis written in terms of the Laguerre polynomial $L_n^{2\nu}(1/x)$. We relegate this treatment along with other cases to section 6.

### 3.4: k = 3

Additional parameter constraints on top of those already given by (8a) are needed to eliminate the three non-tridiagonal terms inside the curly brackets in (9) which are proportional to $x^2$, $x^3$ and $x^4$. The result is the following parameter relations:

$$\alpha = (b-1)\left(\tfrac{a}{2}-1\right) - A_-, \qquad 2\beta = 1-b, \qquad \mu = b\left(\tfrac{a}{2}-1\right) - A_-, \tag{25a}$$

$$A_0 = \left(n+\mu+\tfrac{1}{2}\right)^2 - \tfrac{1}{4}(a-1)^2, \qquad A_+ = 0. \tag{25b}$$

Compatibility of these equations dictates that either $A_0$ is quadratic in $n$ or $A_-$ is linear in $n$. However, as explained in subsection 3.3 above the solution in this case is trivial since it just reverts to the differential equation of the Bessel polynomial $J_n^\mu(x)$. Therefore, this case will be treated in section 6.

## 4. TRA solution class corresponding to (8b)

In this section, we treat the case corresponding to (8b) where $\omega(x) = x^{-1}$ and $2\beta = 1-b$. There are two non-tridiagonal terms inside the curly brackets in (5), which are proportional to $x^2$ and $x^{-1}$ that must be eliminated. The result is the following parameter relations:

$$2\beta = 1-b, \qquad b^2 = 1+4A_1, \qquad A_+ = 0. \tag{26}$$

The first one relates the basis parameter $\beta$ to the differential equation parameters leaving $\alpha$ and $\mu$ to be determined later by physical requirements in section 7. On the other hand, the last two relations are restrictions on the differential equation parameters to obtain a solution. The second relates $b$ to $A_1$ and imposes the reality constraint that $A_1 \geq -\tfrac{1}{4}$. The last relation eliminates the $A_+$ term from the differential equation making the singularity at infinity regular. Therefore, the number of free parameters of the differential equation are reduced from six to four while there are two free basis parameter. With all of these tridiagonal representation constraints, equation (5) transforms into the following form

$$\mathcal{D}\phi_n(x) = G_n x^{\alpha-1} e^{-\beta/x} \left\{ \kappa + x\left[\left(n+\mu+\tfrac{1}{2}\right)^2 - \chi^2\right] + 2\sigma_+ x^2 \frac{d}{dx} \right\} J_n^\mu(x), \tag{27}$$

where $\kappa = \xi + \tfrac{a}{2} + \alpha - 1$, $\chi^2 = \nu^2 + \sigma_+ \sigma_-$ and $\sigma_\pm = -\left(\mu+\tfrac{1}{2}\right) \pm \left(\alpha+\tfrac{a-1}{2}\right)$. Using the recursion relation of the Bessel polynomial (A2) for the second term and the differential relation (A8) for the last term, this equation turns into a form identical to (6) with $\omega(x) = 1/4x$ and giving the following recursion coefficients (for basis normalization $G_n = 1$)

$$u_n = 4\kappa + 2\frac{\mu\chi^2 + 2\sigma_+\left(\mu+\tfrac{1}{2}\right)^2 - (\mu+2\sigma_+)\left(n+\mu+\tfrac{1}{2}\right)^2}{(n+\mu)(n+\mu+1)}, \tag{28a}$$

$$s_n = -\frac{(n+1)\left[\left(n+\mu+\tfrac{3}{2}\right)^2 - \chi^2 - 2\sigma_+(n+2\mu+2)\right]}{(n+\mu+1)\left(n+\mu+\tfrac{3}{2}\right)}, \tag{28b}$$

$$t_n = \frac{(n+2\mu+1)\left[\left(n+\mu+\tfrac{1}{2}\right)^2 - \chi^2 + 2\sigma_+ n\right]}{(n+\mu+1)\left(n+\mu+\tfrac{1}{2}\right)}. \tag{28c}$$

Comparing the recursion relation (7b) constructed from these coefficients to (B7) in Appendix B, we conclude that $P_n$ is a special case of the discrete Hahn polynomial $Q_n^N(k;p,q)$ with

$$p = \alpha + \tfrac{a-1}{2} - 1 \pm \nu, \qquad q = 2\mu + 1 - \left(\alpha + \tfrac{a-1}{2}\right) \mp \nu, \tag{29a}$$

$$N = -\left(\alpha + \tfrac{a-1}{2}\right) \pm \nu, \qquad k = -\kappa. \tag{29b}$$

Therefore, the un-normalized solution of the differential equation (1) for this class reads as follows

$$y_k(x) = \sum_n C_n Q_n^N(k;p,q) \phi_n(x), \tag{30}$$

where $C_n = \left[\dfrac{n + \mu + \tfrac{1}{2}}{\mu + \tfrac{1}{2}} \dfrac{(-n-2\mu)_n}{n!}\right] \prod_{m=0}^{n-1} \dfrac{\left(m+\mu+\tfrac{1}{2}\right)^2 - \chi^2 + 2\sigma_+ m}{\left(m+\mu+\tfrac{3}{2}\right)^2 - \chi^2 - 2\sigma_+(m+2\mu+2)}$.

On the other hand, if we compare the recursion relation (7b) constructed from the coefficients (28) to (B9) but for $(-1)^n H_n^p(z;q,c,d)$, then we get $p = -\kappa$, $q = -\kappa + 2(\mu+1) - (2\alpha + a - 1)$, $c = \kappa + \left(\alpha + \tfrac{a-1}{2}\right) + \nu$, $d = \kappa + \left(\alpha + \tfrac{a-1}{2}\right) - \nu$ and $z = 0$.

## 5. TRA solution class corresponding to (8c)

Finally, we treat the case corresponding to (8c) where $\omega(x) = x^{-2}$ and $\alpha = \mu + 1 - \tfrac{a}{2}$. There are two non-tridiagonal terms inside the curly brackets in (5), which are proportional to $x^2$ and $x^3$, to be eliminated. After a similar discussion that followed (22b) in section 3.3, the result for non-trivial solution requires the following parameter relations:

$$\alpha(n) = -n - \nu + \tfrac{1}{2} - \tfrac{a}{2}, \qquad \mu(n) = -n - \nu - \tfrac{1}{2}, \qquad A_+ = 0, \tag{31}$$

Now, with all of these tridiagonal representation constraints, equation (5) becomes:

$$\mathcal{D}\phi_n(x) = G_n x^{\alpha - 2} e^{-\beta/x} \left\{ A_1 + \beta(\beta + b) + \left[\xi - (1+\tau)\left(n+\nu+\tfrac{1}{2}\right)\right]x + (2\beta + b - 1)x^2 \dfrac{d}{dx}\right\} \bar{J}_n^\nu(x) \tag{32}$$

The comments made below Eq. (24) show that the square bracket term, which is proportional to $x$, destroys the tridiagonal structure. Moreover, the differential property of $\bar{J}_n^\nu(x)$ given by (A13) in addition to (A12) shows that the last term in (32) is not tridiagonal either. Stated differently, since the n-dependent expression in $\phi_n(x)$ for this case contains $L_n^{2\nu}(1/x)$ then the differential relation of the Laguerre polynomial shows that the differential term in (32) destroys the tridiagonal structure because it is the action of $x\tfrac{d}{dx}$ on $L_n^{2\nu}(1/x)$ not $x^2 \tfrac{d}{dx}$ that preserves the tridiagonal structure. Therefore, as done in subsection 3.3 above, we relegate the proper treatment of this case to the following section.

# 6. TRA solution classes in the singular Laguerre basis

The unsuccessful attempt at finding non-trivial TRA solutions for the three cases in subsections 3.3, 3.4 and in section 5 could be resolved by making a proper treatment, which is carried out not in the Bessel basis of Eq. (3) but rather in the singular Laguerre basis to be given below by Eq. (34). In these three cases, $A_+ = 0$ and the Bessel basis parameters $\alpha$ and $\mu$ depend linearly on $n$ as follows

$$\alpha(n) = -n - \nu + \tfrac{1}{2}(1-a), \qquad \mu(n) = -n - \nu - \tfrac{1}{2}, \tag{33}$$

where $\nu = \sqrt{A_0 + \tfrac{1}{4}(a-1)^2}$. The corresponding basis element becomes $G_n x^{-n-\nu-\frac{a}{2}+\frac{1}{2}} e^{-\beta/x} \bar{J}_n^\nu(x)$. Using the definition of $\bar{J}_n^\nu(x)$ given by (A10), we can write this basis element as follows

$$\varphi_n(x) = g_n x^{-\nu-\frac{a}{2}+\frac{1}{2}} e^{-\beta/x} L_n^{2\nu}(1/x), \tag{34}$$

where $g_n = (-1)^n n! G_n$. The well-known recursion relation of the Laguerre polynomial gives [13]

$$\frac{1}{x} L_n^{2\nu}(1/x) = (2n + 2\nu + 1) L_n^{2\nu}(1/x) - (n + 2\nu) L_{n-1}^{2\nu}(1/x) - (n+1) L_{n+1}^{2\nu}(1/x). \tag{35}$$

Moreover, the differential relation of the Laguerre polynomial results in the following

$$x \frac{d}{dx} L_n^{2\nu}(1/x) = -n L_n^{2\nu}(1/x) + (n + 2\nu) L_{n-1}^{2\nu}(1/x). \tag{36}$$

Using the differential equation of the Laguerre polynomial,

$$\left\{ x^2 \frac{d^2}{dx^2} + \left[ 1 + x(1 - 2\nu) \right] + \frac{n}{x} \right\} L_n^{2\nu}(1/x) = 0, \tag{37}$$

we obtain the following action of the second order differential operator $\mathcal{D}$, with $A_+ = 0$ and $\nu^2 = A_0 + \tfrac{1}{4}(a-1)^2$, on the basis elements (34)

$$\mathcal{D}\varphi_n(x) = g_n x^{-\nu-\frac{a}{2}+\frac{1}{2}} e^{-\beta/x} \left\{ (2\beta + b - 1) \frac{d}{dx} + \frac{A_1 + \beta(\beta + b)}{x^2} \right.$$
$$\left. - \frac{1}{x} \left[ -A_- + \tfrac{1}{2} b(a-1) + \beta + \nu(2\beta + b) + n \right] \right\} L_n^{2\nu}(1/x) \tag{38}$$

Taking into account the above properties of the polynomial $L_n^{2\nu}(1/x)$, we arrive at the following three scenarios that result in the tridiagonal structure depicted by Eq. (6):

(1) $2\beta = 1 - b$ and $\omega(x) = x^{-1}$                                                         (39a)
(2) $2\beta = 1 - b$ and $\omega(x) = 1$                                                             (39b)
(3) $2\beta \neq 1 - b$ and $\omega(x) = x^{-1}$.                                                      (39c)

In the following three subsections, we treat these cases separately where we obtain the corresponding recursion relation (7b) and solve it in terms of orthogonal polynomials.

**6.1: Case (39a)**

In addition to the constraints $A_+ = 0$ and $v^2 = A_0 + \frac{1}{4}(a-1)^2$, the tridiagonal requirement dictates that $2\beta = 1 - b$. With all of these parameter constraints, equation (38) transforms into the following:

$$\mathcal{D}\varphi_n(x) = g_n x^{-v-\frac{a}{2}-\frac{1}{2}} e^{-\beta/x} \left\{ \left[ A_1 + \frac{1}{4}(1-b^2) \right] \frac{1}{x} - \left( n + \zeta + \frac{1}{2} \right) \right\} L_n^{2v}(1/x), \tag{40}$$

where $\zeta = -A_- + v + b\left(\frac{a}{2} - 1\right)$. Using the recursion relation (35) for the first term, this equation turns into a form identical to Eq. (6) with the following recursion coefficients

$$u_n = 2 \frac{-2A_- + b(a-2)}{4A_1 - b^2 + 1} - \frac{4A_1 - b^2 - 1}{4A_1 - b^2 + 1}(2n + 2v + 1), \tag{41a}$$

$$s_n = n + 2v + 1, \qquad t_n = n + 1 \tag{41b}$$

Comparing the recursion relation (7b) constructed from these coefficients to (B11) in Appendix B, we conclude that $P_n$ is the Meixner-Pollaczek polynomial $P_n^\lambda(z;\theta)$ with

$$\lambda = v + \tfrac{1}{2}, \qquad \cos\theta = \frac{4A_1 - b^2 - 1}{4A_1 - b^2 + 1}, \qquad z = \frac{2A_- + b(2-a)}{2\sqrt{4A_1 - b^2}}, \tag{42}$$

which requires $4A_1 \geq b^2$. Therefore, the un-normalized solution of the differential equation (1) for this class reads as follows

$$y(x) = \sum_n C_n P_n^\lambda(z;\theta) \varphi_n(x), \tag{43}$$

where $C_n = \left( \prod_{m=0}^{n-1} t_m / s_m \right) = \frac{n!}{(2v+1)_n}$.

**6.2: Case (39b)**

In addition to the constraints $A_+ = 0$ and $v^2 = A_0 + \frac{1}{4}(a-1)^2$, the tridiagonal requirement dictates the following relations among the parameters

$$2\beta = 1 - b, \qquad b^2 = 1 + 4A_1. \tag{44}$$

With all of these parameter constraints, equation (38) transforms into the following:

$$\mathcal{D}\varphi_n(x) = -g_n x^{-v-\frac{a}{2}+\frac{1}{2}} e^{-\beta/x} \left[ \left( n + \zeta + \tfrac{1}{2} \right) \frac{1}{x} \right] L_n^{2v}(1/x). \tag{45}$$

Using the recursion relation (35), this equation turns into a form identical to Eq. (6) with the following recursion coefficients

–11–

$$u_n = -\left(n+\zeta+\tfrac{1}{2}\right)(2n+2\nu+1), \tag{46a}$$

$$s_n = (n+2\nu+1)\left(n+\zeta+\tfrac{3}{2}\right), \qquad t_n = (n+1)\left(n+\zeta+\tfrac{1}{2}\right) \tag{46b}$$

Comparing the recursion relation (7a) constructed from these coefficients to (B5) in Appendix B, we conclude that $Q_n$ (where $f_n = f_0 Q_n$) is the continuous dual Hahn polynomial $S_n^p(z^2;c,d)$ with

$$p = \nu+1, \qquad c = \nu, \qquad d = \zeta-\nu+\tfrac{1}{2}, \qquad z^2 = -\nu^2. \tag{47}$$

which requires that $z$ be imaginary indicating discrete spectrum. Otherwise, $4A_0 < -(a-1)^2$ making $\nu$ pure imaginary. The un-normalized solution of the differential equation (1) for this class reads as follows

$$y(x) = \sum_n S_n^{\nu+1}\left(-\nu^2;\nu,d\right)\varphi_n(x), \tag{48}$$

On the other hand, if we compare the recursion relation (7a) constructed from coefficients (46) to (B3) in Appendix B, we conclude that $Q_n$ (where $f_n = f_0 Q_n$) is the dual Hahn polynomial $R_n^N(z_k^2;p,q)$ with

$$p = \zeta+\tfrac{1}{2}, \qquad q = 2\nu-\zeta+\tfrac{1}{2}, \qquad -N = 2\nu+1, \qquad z_k^2 = \nu^2. \tag{49}$$

The positivity requirement on $\nu$ is in contradiction with the values for $N$. Thus, this alternative solution must be rejected in favor of (48).

### 6.3: Case (39c)

With the constraints $A_+ = 0$ and $\nu^2 = A_0 + \tfrac{1}{4}(a-1)^2$, equation (38) becomes:

$$\mathcal{D}\varphi_n(x) = g_n x^{-\nu-\tfrac{a}{2}-\tfrac{1}{2}} e^{-\beta/x}\left\{\tau x \frac{d}{dx} + \frac{1}{4x}\left[4A_1 - b^2 + (\tau+1)^2\right]\right.$$
$$\left. - \left[\left(n+\zeta+\tfrac{1}{2}\right)+\tau\left(\nu+\tfrac{1}{2}\right)\right]\right\} L_n^{2\nu}(1/x) \tag{50}$$

where $\tau = 2\beta+b-1$. Using the differential relation (36) for the first term and the recursion relation (35) for the second term, this equation turns into a form identical to Eq. (6) with the following recursion coefficients

$$u_n = 2\left[-2A_- + b(a-2)\right] - \left(4A_1 - b^2 + \tau^2 - 1\right)(2n+2\nu+1), \tag{51a}$$

$$s_n = \left[4A_1 - b^2 + (\tau-1)^2\right](n+2\nu+1), \tag{51b}$$

$$t_n = \left[4A_1 - b^2 + (\tau+1)^2\right](n+1). \tag{51c}$$

Note that if $\tau = 0$ this case gives the same recursion relation as that of (39a), which resulted in the Meixner-Pollaczek polynomial $P_n^\lambda(z;\theta)$. Therefore, $\tau$ could be viewed as a deformation parameter. If we divide the recursion relation by $\left(4A_1 - b^2 + \tau^2 + 1\right)$ then we obtain the three-term recursion relation of the polynomial $Y_n^\lambda(z;\theta,\eta)$ shown as Eq. (B16) in Appendix B, where

$$\lambda = v + \tfrac{1}{2}, \quad \cos\theta = \frac{4A_1 - b^2 + \tau^2 - 1}{4A_1 - b^2 + \tau^2 + 1}, \quad \eta = \frac{\tau}{\sqrt{4A_1 - b^2 + \tau^2}}, \quad z = \frac{2A_- + b(2-a)}{2\sqrt{4A_1 - b^2 + \tau^2}}, \quad (52)$$

which requires $4A_1 \geq b^2 - \tau^2$. In fact, the orthogonal polynomial associated with this case depends on the value of $\eta$. Therefore, if $|\eta| > 1$ then it will be the discrete polynomial $Z_n^{v+\frac{1}{2}}(m;\theta,\eta)$ defined at the end of Appendix B. However, if $|\eta| < 1$ then the polynomial will be $Y_n^{v+\frac{1}{2}}(z;\theta,\eta)$ and the un-normalized solution of the differential equation (1) reads as follows

$$y(x) = \sum_n C_n Y_n^{v+\frac{1}{2}}(z;\theta,\eta)\varphi_n(x), \quad (53)$$

where $C_n = \left(\prod_{m=0}^{n-1} t_m/s_m\right) = \dfrac{n!}{(2v+1)_n}\left[\dfrac{1+\eta\sin\theta}{1-\eta\sin\theta}\right]^n$.

## 7. Physical applications

As an application of our solutions to Eq. (1), we present in this section two illustrative examples of quantum systems, identify their potential functions and obtain some of their physical properties. One of these is a novel system in one dimension with a confining exponential potential and the other is a singular isotropic oscillator.

In the atomic units $\hbar = m = 1$, the Schrödinger equation for a potential function $V(r)$ and energy $E$ reads as follows

$$\left[-\frac{1}{2}\frac{d^2}{dr^2} + V(r) - E\right]\psi(r) = 0, \quad (54)$$

where the configuration space coordinate $r$ is either the whole real line $r \in \mathbb{R}$, half of the line $r \geq 0$, or a finite segment of the line $r \in [r_-, r_+]$. If we rewrite $V(r) \mapsto V(r) + \ell(\ell+1)/2r^2$, then Eq. (54) represents the radial Schrödinger equation in three dimensions with spherical symmetry and angular momentum quantum number $\ell$. Now, we make a coordinate transformation $r \mapsto x(\lambda r)$ where $\lambda$ is a positive scale parameter of inverse length dimension such that $x \geq 0$ and $dx/dr = \eta\lambda\left(x^a e^{-b/x}\right)$, where $\eta$ is a real constant parameter. This transformation is compatible with the Bessel polynomial weight function (A3) and maps the above Schrödinger equation into the following second order differential equation in terms of the dimensionless variable $x$

$$x^{2a-2}e^{-2b/x}\left\{x^2\frac{d^2}{dx^2}+(ax+b)\frac{d}{dx}+x^{2-2a}e^{2b/x}\left[\varepsilon-U(r)\right]\right\}\psi(r)=0, \tag{55}$$

where $\varepsilon=\left(2/\eta^2\lambda^2\right)E$ and $U(r)=\left(2/\eta^2\lambda^2\right)V(r)$. Since the factor to the left, $x^{2a-2}e^{-2b/x}$, is a positive definite function for $x>0$ then the solution of this equation is the same as that of the Bessel-type equation (1) provided that we can write

$$\varepsilon-U(r)=x^{2a-2}e^{-2b/x}\left[A_+x+\frac{A_-}{x}+\frac{A_1}{x^2}-A_0\right]. \tag{56}$$

To obtain a non-zero energy solution, the right-hand side of this equation must contain a constant to be identified with the energy on the left. To achieve that, we must take $b=0$ and choose $a$ from one of the four possibilities $a=\{\tfrac{1}{2},1,\tfrac{3}{2},2\}$ leading to the four scenarios shown in Table 1. As illustration, we give two examples: one that belongs to the class of solutions in subsection 3.1 and another that belongs to the class of solutions in subsection 6.3.

**Table 1**: List of quantum mechanical systems whose Schrödinger equation could be transformed into the Bessel-type equation (1). The configuration space coordinate $r$ is related to the differential equation variable $x$ by $dx/dr=\eta\lambda\left(x^a e^{-b/x}\right)$ with $b=0$. The effective potential function is $V(r)+\frac{\ell(\ell+1)+\Lambda}{2r^2}$.

| $a$ | $x(r)$ | $\eta$ | $E$ | $V(r)$ | $\ell(\ell+1)+\Lambda$ |
|---|---|---|---|---|---|
| $\tfrac{1}{2}$ | $(\lambda r)^2$, $r\geq 0$ | 2 | $2\lambda^2 A_+$ | $-\dfrac{2A_-/\lambda^2}{r^4}-\dfrac{2A_1/\lambda^4}{r^6}$ | $4A_0$ |
| 1 | $e^{\lambda r}$, $-\infty<r<+\infty$ | 1 | $-\lambda^2 A_0/2$ | $-\dfrac{\lambda^2}{2}\left(A_+e^{\lambda r}+A_-e^{-\lambda r}+A_1 e^{-2\lambda r}\right)$ | 0 |
| $\tfrac{3}{2}$ | $(\lambda r)^{-2}$, $r\geq 0$ | $-2$ | $2\lambda^2 A_-$ | $-\dfrac{2A_+/\lambda^2}{r^4}-\left(2\lambda^4 A_1\right)r^2$ | $4A_0$ |
| 2 | $(\lambda r)^{-1}$, $r\geq 0$ | $-1$ | $\lambda^2 A_1/2$ | $-\dfrac{\lambda A_-/2}{r}-\dfrac{A_+/2\lambda}{r^3}$ | $A_0$ |

For the first example, we take $(a,b)=(1,0)$ giving the coordinate transformation $x(\lambda r)=e^{\lambda r}$, where $-\infty<r<+\infty$. According to Eq. (10), the problem parameters are

$$\alpha=\tfrac{1}{2}-A_-, \quad \beta=\tfrac{1}{2}, \quad \mu=-A_-, \quad A_1=-\tfrac{1}{4}, \tag{57}$$

which requires that $A_-\geq N+\tfrac{1}{2}$. Moreover, the potential function and energy are obtained from Eq. (56) as follows

–14–

$$V(r) = \frac{\lambda^2}{2}\left(\tfrac{1}{4}e^{-2\lambda r} - A_-e^{-\lambda r} - A_+e^{\lambda r}\right), \qquad E = -\frac{\lambda^2}{2}A_0. \qquad (58)$$

Physically, we require that $A_+ \leq 0$ otherwise the particle escapes to $+\infty$. If $A_+ = 0$ then the problem becomes the usual 1D Morse potential [7,8]. On the other hand, if $A_+ < 0$ then the problem becomes that of a one-dimensional exponentially confining potential well. This problem has been treated recently by one of the authors [16]. The wavefunction is given by Eq. (15) as

$$\psi(x) = \sum_n C_n B_n^{-A_-}\left(4v^2/A_+ ; 4/A_+\right)\phi_n(x), \qquad (59)$$

where $\phi_n(x) = G_n x^{\frac{1}{2}-A_-} e^{-1/2x} J_n^{-A_-}(x)$ and $C_n = \left(\prod_{m=0}^{n-1} t_m/s_m\right) = \left[\frac{n - A_- + \frac{1}{2}}{\frac{1}{2} - A_-}\right]\frac{(-n + 2A_-)_n}{n!}$. Since the analytic properties of the new polynomial $B_n^\mu(z;\gamma)$ are not known yet, numerical means were used in [16] to obtain the physical properties of the system such as the energy spectrum. Interested readers can find more information about this system in [16].

The second example belongs to subsection 6.3 where we take $(a,b) = \left(\tfrac{3}{2}, 0\right)$ giving the coordinate transformation $x(\lambda r) = (\lambda r)^{-2}$, where $r \geq 0$. The problem parameters are

$$A_+ = 0, \qquad 2\beta = \tau + 1, \qquad v^2 = A_0 + \tfrac{1}{16}, \qquad (60)$$

which requires $16A_0 \geq -1$. Moreover, the effective radial potential function and energy are obtained from Eq. (56) as follows

$$\frac{\ell(\ell+1)}{2r^2} + \frac{\Lambda}{2r^2} - \left(2\lambda^4 A_1\right)r^2, \qquad E = 2\lambda^2 A_-, \qquad (61)$$

where $\Lambda = 4A_0 - \ell(\ell+1) = 4v^2 - \left(\ell + \tfrac{1}{2}\right)^2$. Physically, we require that $A_1 \leq 0$ otherwise the particle escapes to $+\infty$. Moreover, to avoid the quantum anomalies associated with the singular $r^{-2}$ potential [17,18] ("fall to the center" problem [19]) we require that $\Lambda \geq -\left(\ell + \tfrac{1}{2}\right)^2$, which is the condition $16A_0 \geq -1$ that has just been stated above. Note that if $\Lambda = 0$, then the problem becomes the conventional isotropic oscillator [7,8]. On the other hand, if $\Lambda \geq -\left(\ell + \tfrac{1}{2}\right)^2$ then we can write the wavefunction as given by (53) that reads

$$\psi(r) = \sum_n C_n Y_n^{v+\frac{1}{2}}(z;\theta,\eta)\varphi_n(r), \qquad (62)$$

where $\varphi_n(r) = g_n(\lambda r)^{2v+\frac{1}{2}} e^{-(\tau+1)\lambda^2 r^2/2} L_n^{2v}\left(\lambda^2 r^2\right)$. The wavefunction parameters are obtained from Eq. (52) as follows

$$\cos\theta = \frac{4A_1 + \tau^2 - 1}{4A_1 + \tau^2 + 1}, \qquad \eta = \frac{\tau}{\sqrt{4A_1 + \tau^2}}, \qquad z = \frac{A_-}{\sqrt{4A_1 + \tau^2}}. \qquad (63)$$

–15–

The energy spectrum is obtained from the spectrum formula of the polynomial $Y_n^{v+\frac{1}{2}}(z;\theta,\eta)$, which in turn is obtained from the spectrum formula of the Meixner-Pollaczek polynomial $P_n^{v+\frac{1}{2}}(\tilde{z};\tilde{\theta})$ together with relation (B19). Now, the spectrum formula for $P_n^{v+\frac{1}{2}}(\tilde{z};\tilde{\theta})$ is $\tilde{z}_k^2 = -\left(k+v+\tfrac{1}{2}\right)^2$ for $k = 0,1,2,..$ (see Appendix A.1 of [14]). Therefore, using (B19) and $z = \tilde{z}\sqrt{1-\eta^2}$ we obtain the following energy spectrum formula for this singular isotropic oscillator

$$E_k = 4\lambda^2 \sqrt{-A_1}\left[k+\tfrac{1}{2}+\tfrac{1}{2}\sqrt{\Lambda+\left(\ell+\tfrac{1}{2}\right)^2}\right]. \tag{64}$$

Finally, we conclude the work with a summary of our findings and discussion of some relevant points.

## 8. Conclusion and discussion

In this paper, we continued our program of using the Tridiagonal Representation Approach to solve physically significant linear ordinary differential equation of the second order. These equation are written in the following generic form

$$\left[p(x)\frac{d^2}{dx^2}+q(x)\frac{d}{dx}+F(x)\right]y(x)=0, \tag{65}$$

where $p(x)$ and $q(x)$ are polynomials of at most degree three. The function $F(x)$, on the other hand, is highly non-trivial having the following form

$$F(x)=\frac{C_r}{r-x}+\frac{C_-}{1-x}+\frac{C_+}{x}+C_1 x+C_0, \tag{66}$$

where $\{r,C_i\}$ are real parameters. Table 2 lists the five differential equations that have been treated by us so far. All solution classes of these equations are written as convergent series (finite or infinite) of square integrable functions. The expansion coefficients of the series are orthogonal polynomials in the parameters space of the corresponding differential equation. The properties of the solutions are obtained from the properties of these polynomials. However, in the course of our endeavor we have identified many such polynomials as deformation of known orthogonal polynomials that are very rich in structure. Examples, include the deformed Meixner-Pollaczek polynomial (B16), the deformed Bessel polynomials (B1) as well as others [20]. As a physical application, one could transform the wave equation to take the form given by Eq. (65). Consequently, the wavefunction could be written in terms of the established solution $y(x)$ and the corresponding interaction potential and energy could be deduced from the function $F(x)$ and the type of transformation of the wavefunction and configuration space coordinate. Moreover, the physical properties of the system (e.g., energy spectrum of the bound states, phase shift of the scattering states, density of states, etc.) could be obtained from the properties of these orthogonal polynomials (e.g., weight function, generating function, asymptotics, zeros, etc.).

Table 2: List of the five differential equations that have been treated by us so far using the Tridiagonal Representation Approach. The general form of these equations is given by Eq. (65).

| Equation Type [Reference] | $p(x)$ | $q(x)$ | $F(x)$ |
|---|---|---|---|
| Confluent Hypergeometric (Laguerre-type) [9] | $x$ | Linear | $\frac{C_+}{x}+C_1 x+C_0$ |
| Hypergeometric (Jacobi-type) [9] | $x(1-x)$ | Linear | $\frac{C_-}{1-x}+\frac{C_+}{x}+C_1 x+C_0$ |
| Bessel-type [this work] | $x^2$ | Linear | $\frac{C_2}{x^2}+\frac{C_+}{x}+C_1 x+C_0$ |
| Heun-type [10] | $x(1-x)(r-x)$ | Quadratic | $\frac{C_r}{r-x}+\frac{C_-}{1-x}+\frac{C_+}{x}+C_1 x+C_0$ |
| Generalized Heun-type [11] | $x(1-x)(r-x)$ | Cubic | $\frac{C_r}{r-x}+\frac{C_-}{1-x}+\frac{C_+}{x}+C_1 x+C_0$ |

## Acknowledgments

This work is supported by the Saudi Center for Theoretical Physics (SCTP). The support by King Fahd University of Petroleum and Minerals (KFUPM) is acknowledged by both authors.

## Appendix A: The Bessel polynomial on the real line

The Bessel polynomial on the real line is defined in terms of the hypergeometric or confluent hypergeometric functions as follows (see section 9.13 of the book by Koekoek et. al [12]):

$$J_n^\mu(x) = {}_2F_0\left(\begin{array}{c}-n,n+2\mu+1\\ -\end{array}\bigg| -x\right) = (n+2\mu+1)_n \, x^n \, {}_1F_1\left(\begin{array}{c}-n\\ -2(n+\mu)\end{array}\bigg| 1/x\right), \tag{A1}$$

where $x \geq 0$, $n = 0,1,2,..,N$ and $N$ is a non-negative integer. The real parameter $\mu$ is negative such that $\mu < -N - \frac{1}{2}$. The Pochhammer symbol $(a)_n$ (a.k.a. shifted factorial) is defined as $(a)_n = a(a+1)(a+2)...(a+n-1) = \frac{\Gamma(n+a)}{\Gamma(a)}$. The Bessel polynomial could also be written in terms of the associated Laguerre polynomial as: $J_n^\mu(x) = n!(-x)^n L_n^{-(2n+2\mu+1)}(1/x)$. The three-term recursion relation reads as follows:

$$\begin{aligned}2x J_n^\mu(x) =& \frac{-\mu}{(n+\mu)(n+\mu+1)} J_n^\mu(x) \\ &- \frac{n}{(n+\mu)(2n+2\mu+1)} J_{n-1}^\mu(x) + \frac{n+2\mu+1}{(n+\mu+1)(2n+2\mu+1)} J_{n+1}^\mu(x)\end{aligned} \tag{A2}$$

Note that the constraints on $\mu$ and on the maximum polynomial degree make this recursion definite (i.e., the signs of the two recursion coefficients multiplying $J^{\mu}_{n\pm1}(x)$ are the same). Otherwise, these polynomials could not be defined on the real line but on the unit circle in the complex plane. The orthogonality relation reads as follows

$$\int_0^\infty x^{2\mu} e^{-1/x} J^{\mu}_n(x) J^{\mu}_m(x)\, dx = -\frac{n!\Gamma(-n-2\mu)}{2n+2\mu+1}\delta_{nm}. \tag{A3}$$

The differential equation is

$$\left\{ x^2 \frac{d^2}{dx^2} + [1+2x(\mu+1)]\frac{d}{dx} - n(n+2\mu+1) \right\} J^{\mu}_n(x) = 0. \tag{A4}$$

The forward and backward shift differential relations read as follows

$$\frac{d}{dx} J^{\mu}_n(x) = n(n+2\mu+1) J^{\mu+1}_{n-1}(x). \tag{A5}$$

$$x^2 \frac{d}{dx} J^{\mu}_n(x) = -(2\mu x+1) J^{\mu}_n(x) + J^{\mu-1}_{n+1}(x). \tag{A6}$$

We can write $J^{\mu-1}_{n+1}(x)$ in terms of $J^{\mu}_n(x)$ and $J^{\mu}_{n\pm1}(x)$ as follows

$$2 J^{\mu-1}_{n+1}(x) = \frac{(n+1)(n+2\mu)}{(n+\mu)(n+\mu+1)} J^{\mu}_n(x)$$
$$+ \frac{n(n+1)}{(n+\mu)(2n+2\mu+1)} J^{\mu}_{n-1}(x) + \frac{(n+2\mu)(n+2\mu+1)}{(n+\mu+1)(2n+2\mu+1)} J^{\mu}_{n+1}(x) \tag{A7}$$

Using this identity and the recursion relation (A2), we can rewrite the backward shift differential relation as follows

$$2x^2 \frac{d}{dx} J^{\mu}_n(x) = n(n+2\mu+1) \times$$
$$\left[ -\frac{J^{\mu}_n(x)}{(n+\mu)(n+\mu+1)} + \frac{J^{\mu}_{n-1}(x)}{(n+\mu)(2n+2\mu+1)} + \frac{J^{\mu}_{n+1}(x)}{(n+\mu+1)(2n+2\mu+1)} \right] \tag{A8}$$

The generating function is

$$\sum_{n=0}^{\infty} J^{\mu}_n(x) \frac{t^n}{n!} = \frac{2^{2\mu}}{\sqrt{1-4xt}} \left(1+\sqrt{1-4xt}\right)^{-2\mu} \exp\left[2t/(1+\sqrt{1-4xt})\right]. \tag{A9}$$

In subsections 3.3, 3.4 and in section 5 of the paper, the polynomial parameter $\mu$ depends on $n$ as $\mu(n) = -n-\nu-\tfrac{1}{2}$. In this case, we can rewrite the polynomial as follows

$$J^{\mu(n)}_n(x) := \bar{J}^{\nu}_n(x) = {}_2F_0\!\left(\begin{array}{c}-n,-n-2\nu\\ -\end{array}\Big| -x\right)$$
$$= (2\nu+1)_n (-x)^n \, {}_1F_1\!\left(\begin{array}{c}-n\\ 2\nu+1\end{array}\Big| 1/x\right) = n!(-x)^n L^{2\nu}_n(1/x) \tag{A10}$$

–18–

The corresponding recursion relation (A2) changes to

$$2\nu(4\nu^2 - 1)x\bar{J}_n^\nu(x) = 2\nu(2n + 2\nu + 1)\bar{J}_n^\nu(x) - (2\nu - 1)n\bar{J}_{n-1}^{\nu+1}(x) - (2\nu + 1)(n + 2\nu)\bar{J}_{n+1}^{\nu-1}(x). \quad (A11)$$

This is not a typical three-term recursion relation for $\bar{J}_n^\nu(x)$. Using the properties of the Laguerre polynomials, we can show that $\bar{J}_{n\pm 1}^{\nu\mp 1}(x)$ is a combination of $\bar{J}_n^\nu(x)$ and $\bar{J}_{n\pm 1}^\nu(x)$ but with $x$ dependent factors as follows

$$\bar{J}_{n+1}^{\nu-1}(x) = \bar{J}_{n+1}^\nu(x) + 2(n+1)x\bar{J}_n^\nu(x) + n(n+1)x^2\bar{J}_{n-1}^\nu(x), \quad (A12a)$$

$$\bar{J}_{n-1}^{\nu+1}(x) = \bar{J}_{n+1}^\nu(x) + 2(n + 2\nu + 1)x\bar{J}_n^\nu(x) + (n + 2\nu)(n + 2\nu + 1)x^2\bar{J}_{n-1}^\nu(x). \quad (A12b)$$

The backward shift differential relation (A8) changes as follows

$$2\nu(4\nu^2 - 1)x^2 \frac{d}{dx}\bar{J}_n^\nu(x) = n(n + 2\nu)\left[4\nu\bar{J}_n^\nu(x) - (2\nu - 1)\bar{J}_{n-1}^{\nu+1}(x) - (2\nu + 1)\bar{J}_{n+1}^{\nu-1}(x)\right]. \quad (A13)$$

The new differential equation replacing (A4) is

$$\left\{x^2 \frac{d^2}{dx^2} + \left[1 - x(2n + 2\nu - 1)\right]\frac{d}{dx} + n(n + 2\nu)\right\}\bar{J}_n^\nu(x) = 0. \quad (A14)$$

## Appendix B: Orthogonal polynomials as expansion coefficients of the solution series

Symbols used in this Appendix are local and should not be confused with those in the rest of the paper. In [16], we introduced a new orthogonal polynomial, which we referred to as $B_n^\mu(x;\gamma)$. It is in fact a deformation or generalization of the Bessel polynomial with $J_n^\mu(x) = B_n^\mu(4x;0)$. It satisfies the following three-term recursion relation

$$xB_n^\mu(x;\gamma) = \left[\frac{-2\mu}{(n+\mu)(n+\mu+1)} + \gamma\left(n + \mu + \tfrac{1}{2}\right)^2\right]B_n^\mu(x;\gamma)$$
$$- \frac{n}{(n+\mu)(n+\mu+\tfrac{1}{2})}B_{n-1}^\mu(x;\gamma) + \frac{n+2\mu+1}{(n+\mu+1)(n+\mu+\tfrac{3}{2})}B_{n+1}^\mu(x;\gamma) \quad (B1)$$

where $\gamma$ is the deformation parameter, $n = 0,1,2,..,N$ and $N$ is a non-negative integer. The real parameter $\mu$ is negative such that $\mu < -N - \tfrac{1}{2}$. This type of deformation/generalization is, in fact, common to all polynomials in the Askey scheme of hypergeometric orthogonal polynomials (see, for example section 5 of [20] and section IV of [21]). Unfortunately, the analytic properties of these polynomials (e.g., weight functions, generating functions, asymptotics, zeros, etc.) are not yet known. This is an open problem in orthogonal polynomials along with other similar ones. For discussions about these open problems, one may consult References [20,22,23]. Nonetheless, these polynomials can be written explicitly for any degree (albeit not in closed form) using the recursion (B1) and starting from the initial seed $B_0^\mu(x;\gamma) = 1$ and setting $B_{-1}^\mu(x;\gamma) \equiv 0$. Had these properties been known, we would have been able to extract all physical information of the corresponding system (e.g., energy spectrum, scattering phase shift, density of states, etc.). That is why we had to resort to numerical means to obtain the physical results as noted in section 7.

In all cases except two (*cf.* subsections 3.1 and 6.2) the recursion coefficients $s_n$ and $t_n$ are either quadratic polynomials in $n$ or a quotient of cubic over quadratic polynomials in $n$. The former situation corresponds to the dual Hahn or continuous dual Hahn polynomials whereas the latter corresponds to the Hahn or continuous Hahn polynomials (see sections 9.3 to 9.6 in Ref. [12]). For ease of reference, we define these polynomials and write their three-term recursion relations. Other properties could be found elsewhere (e.g., in Ref. [12]).

The three-parameter discrete dual Hahn polynomial is defined as (see section 9.6 of [12])

$$R_n^N(z_m^2; p, q) = {}_3F_2\left(\begin{matrix}-n,-m,m+p+q+1\\p+1,-N\end{matrix}\bigg|1\right), \quad (B2)$$

where $n, m = 0, 1, 2, .., N$ and $z_m^2 = \left(m + \frac{p+q+1}{2}\right)^2$. The parameters $p$ and $q$ are either greater than $-1$ or less than $-N$. The three-term recursion relation is as follows

$$z_m^2 R_n^N(z_m^2; p, q) = \left[(N-n)(n+p+1) + n(N-n+q+1) + \tfrac{1}{4}(p+q+1)^2\right] R_n^N(z_m^2; p, q) \\ - n(N-n+q+1) R_{n-1}^N(z_m^2; p, q) - (N-n)(n+p+1) R_{n+1}^N(z_m^2; p, q) \quad (B3)$$

The three-parameter continuous dual Hahn polynomial is defined as (see section 9.3 of [12])

$$S_n^p(x^2; c, d) = {}_3F_2\left(\begin{matrix}-n, p+ix, p-ix\\p+c, p+d\end{matrix}\bigg|1\right), \quad (B4)$$

where either $\mathrm{Re}(p, c, d) > 0$ or $p < 0$ with $p+c > 0$ and $p+d > 0$. Non-real parameters are in conjugate pairs. It satisfies the following three-term recursion relation

$$x^2 S_n^p(x^2; c, d) = \left[n(n+c+d-1) + (n+p+c)(n+p+d) - p^2\right] S_n^p(x^2; c, d) \\ - n(n+c+d-1) S_{n-1}^p(x^2; c, d) - (n+p+c)(n+p+d) S_{n+1}^p(x^2; c, d) \quad (B5)$$

The three-parameter discrete Hahn polynomial is defined as (see section 9.5 of [12])

$$Q_n^N(m; p, q) = {}_3F_2\left(\begin{matrix}-n,-m,n+p+q+1\\p+1,-N\end{matrix}\bigg|1\right), \quad (B6)$$

where $n, m = 0, 1, 2, .., N$. It satisfies the following three-term recursion relation

$$m Q_n^N(m; p, q) = \left[\frac{n(n+q)(n+p+q+N+1)}{(2n+p+q)(2n+p+q+1)} + \frac{(N-n)(n+p+1)(n+p+q+1)}{(2n+p+q+1)(2n+p+q+2)}\right] Q_n^N(m; p, q) \\ - \frac{n(n+q)(n+p+q+N+1)}{(2n+p+q)(2n+p+q+1)} Q_{n-1}^N(m; p, q) - \frac{(N-n)(n+p+1)(n+p+q+1)}{(2n+p+q+1)(2n+p+q+2)} Q_{n+1}^N(m; p, q) \quad (B7)$$

The four-parameter continuous Hahn polynomial is defined as (see section 9.4 of [12])

$$H_n^p(x;q,c,d) = {}_3F_2\left(\begin{matrix}-n,,n+p+q+c+d-1,p+ix\\p+c,p+d\end{matrix}\bigg|1\right), \tag{B8}$$

where $\text{Re}(p,q,c,d) > 0$, $c = \bar{p}$ and $d = \bar{q}$. It satisfies the following three-term recursion relation

$$(p+ix)H_n^p(x;q,c,d) = $$
$$\left[\frac{(n+p+c)(n+p+d)(n+p+q+c+d-1)}{(2n+p+q+c+d-1)(2n+p+q+c+d)} - \frac{n(n+q+c-1)(n+q+d-1)}{(2n+p+q+c+d-2)(2n+p+q+c+d-1)}\right]H_n^p \tag{B9}$$
$$+\frac{n(n+q+c-1)(n+q+d-1)H_{n-1}^p}{(2n+p+q+c+d-2)(2n+p+q+c+d-1)} - \frac{(n+p+c)(n+p+d)(n+p+q+c+d-1)H_{n+1}^p}{(2n+p+q+c+d-1)(2n+p+q+c+d)}$$

For the rest of this Appendix, we consider orthogonal polynomials that are relevant to some of the cases treated in section 6. We start with the two-parameter Meixner-Pollaczek polynomial, which is defined as (see section 9.7 of [12])

$$P_n^\lambda(x;\theta) = \frac{(2\lambda)_n}{n!}e^{in\theta}\,{}_2F_1\left(\begin{matrix}-n,\lambda+ix\\2\lambda\end{matrix}\bigg|1-e^{-2i\theta}\right), \tag{B10}$$

where $\lambda > 0$ and $0 < \theta < \pi$. It satisfies the following three-term recursion relation

$$2x(\sin\theta)P_n^\lambda(x;\theta) = -2(n+\lambda)(\cos\theta)P_n^\lambda(x;\theta)$$
$$+(n+2\lambda-1)P_{n-1}^\lambda(x;\theta)+(n+1)P_{n+1}^\lambda(x;\theta) \tag{B11}$$

The generating function is

$$\sum_{n=0}^\infty P_n^\lambda(x;\theta)t^n = \left(1-te^{i\theta}\right)^{-\lambda+ix}\left(1-te^{-i\theta}\right)^{-\lambda-ix}. \tag{B12}$$

The discrete version of this polynomial is the Meixner polynomial, which is defined here as (see section 9.10 of Ref. [12] for the conventional definition)

$$M_n^\lambda(m;\theta) = \frac{(2\lambda)_n}{n!}e^{-n\theta}\,{}_2F_1\left(\begin{matrix}-n,-m\\2\lambda\end{matrix}\bigg|1-e^{2\theta}\right), \tag{B13}$$

where $\lambda > 0$, $\theta > 0$ and $n,m = 0,1,2,\ldots$. It is obtained from (B10) by the map $\theta \mapsto i\theta$ and $x \mapsto i(\lambda+m)$. Consequently, under this map (B11) changes to the following three-term recursion relation

$$2m(\sinh\theta)M_n^\lambda(m;\theta) = 2\left[(n+\lambda)\cosh\theta - \lambda\sinh\theta\right]M_n^\lambda(m;\theta)$$
$$-(n+2\lambda-1)M_{n-1}^\lambda(m;\theta)-(n+1)M_{n+1}^\lambda(m;\theta) \tag{B14}$$

Moreover, the generating function (B12) becomes

$$\sum_{n=0}^\infty M_n^\lambda(m;\theta)t^n = \left(1-te^\theta\right)^m\left(1-te^{-\theta}\right)^{-m-2\lambda}. \tag{B15}$$

We also define a deformed version of the Meixner-Pollaczek polynomial, which we designate as $Y_n^\lambda(x;\theta,\eta)$, by its three-term recursion relation that reads as follows

$$2x(\sin\theta)Y_n^\lambda(x;\theta,\eta) = -2(n+\lambda)(\cos\theta)Y_n^\lambda(x;\theta,\eta)$$
$$+(n+2\lambda-1)(1-\eta\sin\theta)Y_{n-1}^\lambda(x;\theta,\eta)+(n+1)(1+\eta\sin\theta)Y_{n+1}^\lambda(x;\theta,\eta) \quad \text{(B16)}$$

where $\eta$ is the deformation parameter. Obviously, $P_n^\lambda(x;\theta) = Y_n^\lambda(x;\theta,0)$. Using this recursion relation, we can derive the generating function as [24]

$$\sum_{n=0}^{\infty} Y_n^\lambda(x;\theta,\eta)t^n = (1-\alpha t)^{-A}(1-\beta t)^{-B}, \quad \text{(B17)}$$

where the parameters $\{\alpha, \beta, A, B\}$ are

$$\alpha = \frac{\cos\theta \pm \sqrt{\eta^2-1}\sin\theta}{1+\eta\sin\theta}, \qquad \beta = \frac{\cos\theta \mp \sqrt{\eta^2-1}\sin\theta}{1+\eta\sin\theta}, \quad \text{(B18a)}$$

$$A = \lambda \pm \frac{x}{\sqrt{\eta^2-1}}, \qquad B = \lambda \mp \frac{x}{\sqrt{\eta^2-1}}. \quad \text{(B18b)}$$

If $|\eta|<1$, then the generating function (B17) shows that

$$Y_n^\lambda(x;\theta,\eta) = \left[\frac{1-\eta\sin\theta}{1+\eta\sin\theta}\right]^{n/2} P_n^\lambda(y;\varphi), \quad \text{(B19)}$$

where $y = \frac{x}{\sqrt{1-\eta^2}}$ and $\cos\varphi = \frac{\cos\theta}{\sqrt{1-\eta^2\sin^2\theta}}$. Moreover, the orthogonality is obtained using that of $P_n^\lambda(x;\theta)$ as follows

$$\frac{(2\sin\varphi)^{2\lambda}}{2\pi\sqrt{1-\eta^2}} \int_{-\infty}^{+\infty} e^{(2\varphi-\pi)y} |\Gamma(\lambda+iy)|^2 Y_n^\lambda(x;\theta,\eta)Y_n^\lambda(x;\theta,\eta)\,dx =$$
$$\left(\frac{1-\eta\sin\theta}{1+\eta\sin\theta}\right)^n \frac{\Gamma(n+2\lambda)}{n!}\delta_{n,m} \quad \text{(B20)}$$

If, on the other hand, $|\eta|>1$ then we perform the map $\theta \mapsto i\theta$, $\eta \mapsto -i\eta$ and $x \mapsto i(\lambda+m)$ in (B16) to obtain the following three-term recursion relation for $Z_n^\lambda(x;\theta,\eta)$, which is the discrete version of $Y_n^\lambda(x;\theta,\eta)$,

$$2m(\sinh\theta)Z_n^\lambda(m;\theta,\eta) = 2[(n+\lambda)\cosh\theta - \lambda\sinh\theta]Z_n^\lambda(m;\theta,\eta)$$
$$-(n+2\lambda-1)(1-\eta\sinh\theta)Z_{n-1}^\lambda(m;\theta,\eta)-(n+1)(1+\eta\sinh\theta)Z_{n+1}^\lambda(m;\theta,\eta) \quad \text{(B21)}$$

Using this recursion relation, one can derive the generating function as

$$\sum_{n=0}^{\infty} Z_n^\lambda(m;\theta,\eta)t^n = (1-te^\varphi)^{\tilde{m}}(1-te^{-\varphi})^{-\tilde{m}-2\lambda}, \quad \text{(B22)}$$

where $\cosh\varphi = \dfrac{\cosh\theta}{\sqrt{1-\eta^2\sinh^2\theta}}$ and $\tilde{m} = \dfrac{m+\lambda}{\sqrt{1+\eta^2}} - \lambda$. Finally, one can show that $Z_n^\lambda(m;\theta,\eta) = \left[\dfrac{1-\eta\sinh\theta}{1+\eta\sinh\theta}\right]^{n/2} M_n^\lambda(\tilde{m};\varphi)$.

## Data Availability Statement:

Data sharing is not applicable to this article as no new data were created or analyzed in this study.

## References:

[1] A. D. Alhaidari and H. Bahlouli, "Tridiagonal Representation Approach in Quantum Mechanics", Phys. Scripta **94**, 125206 (2019)

[2] E. J. Heller and H. A. Yamani, *New $L^2$ Approach to Quantum Scattering: Theory*, Phys. Rev. A **9**, 1201 (1974)

[3] H. A. Yamani and L. Fishman, *J-Matrix Method: Extensions to Arbitrary Angular Momentum and to Coulomb Scattering*, J. Math. Phys. **16**, 410 (1975)

[4] M. E. H. Ismail and E. Koelink, *The J-Matrix method*, Adv. Appl. Math. **46** (2011) 379

[5] M. E. H. Ismail and E. Koelink, *Spectral properties of operators using tridiagonalization*, Analysis and Applications **10** (2012) 327

[6] V. X. Genest, M. E. H. Ismail, L. Vinet, and A. Zhedanov, *Tridiagonalization of the hypergeometric operator and the Racah-Wilson algebra*, Proc. Amer. Math. Soc. **144** (2016) 4441

[7] De R., Dutt R., Sukhatme U., *Mapping of shape invariant potentials under point canonical transformations*, J. Phys. A **25** (1992) L843

[8] F. Cooper, A. Khare and U. Sukhatme, *Supersymmetry in Quantum Mechanics* (World Scientific, Singapore, 2004).

[9] A. D. Alhaidari, *Series solutions of Laguerre- and Jacobi-type differential equations in terms of orthogonal polynomials and physical applications*, J. Math. Phys. **59** (2018) 063508

[10] A. D. Alhaidari, *Series solutions of Heun-type equation in terms of orthogonal polynomials*, J. Math. Phys. **59** (2018) 113507

[11] A. D. Alhaidari, *Series solution of a ten-parameter second order differential equation with three regular singularities and one irregular singularity*, Theor. Math. Phys. **202** (2020) 17 [Russian: TMP **202** (2020) 20]

[12] R. Koekoek, P. A. Lesky and R. F. Swarttouw: *Hypergeometric Orthogonal Polynomials and Their q-Analogues* (Springer, Heidelberg 2010)

[13] T. S. Chihara, *An introduction to orthogonal polynomials*, Mathematics and its Applications, Vol. 13, (Gordon and Breach Science Publishers, New York - London - Paris, 1978)

[14] A. D. Alhaidari and M. E. H. Ismail, *Quantum mechanics without potential function*, J. Math. Phys. **56** (2015) 072107

[15] A. D. Alhaidari, *Representation of the quantum mechanical wavefunction by orthogonal polynomials in the energy and physical parameters*, Commun. Theor. Phys. **72** (2020) 015104

[16] A. D. Alhaidari, *Exponentially confining potential well*, arXiv:2005.09080 [quant-ph], submitted

[17] S. A. Coon and B. R. Holstein, *Anomalies in quantum mechanics: The $1/r^2$ potential*, Am. J. Phys. **70** (2002) 513

[18] A. M. Essin and D. J. Griffiths, *Quantum mechanics of the $1/x^2$ potential*, Am. J. Phys. **74** (2006) 109

[19] L. D. Landau and E. M. Lifshitz, *Quantum Mechanics,* Course of Theoretical Physics Vol. 3, 3rd ed. (Pergamon Press, Oxford, 1977) pp. 114-117

[20] A. D. Alhaidari, *Open problem in orthogonal polynomials*, Rep. Math. Phys. **84** (2019) 393

[21] F. A. Grünbaum, L. Vinet, and A. Zhedanov, *Tridiagonalization and the Heun equation*, J. Math. Phys. **58** (2017) 031703

[22] A. D. Alhaidari, *Orthogonal polynomials derived from the tridiagonal representation approach*, J. Math. Phys. **59** (2018) 013503

[23] W. Van Assche, *Solution of an open problem about two families of orthogonal polynomials*, SIGMA **15** (2019) 005

[24] Private communication with Mourad E. H. Ismail.